# Range of magic constant on Hexagonal Tortoise Problem

Donghwi Park
Korea University, Seoul
twitterion@korea.ac.kr

## 1. Abstract

Hexagonal tortoise problem (HTP), also known as Jisuguimundo or Jisugwimundo, is a magic square variety which was invented by medieval Korean Mathematician and minister Suk-Jung Choi (1646-1715).[1] Choi showed pattern 30 vertices 3 by 3 diagonal shape which has 93 as its magic constant. Unlike magic square, vertices in Jisugwimundo counted one times, twice or three times. This change makes magic constant of Hexagonal Tortoise Problem could be vary. We consider a range of hexagonal sums in various Jisugwimundo. In this paper, we decomposed vertices on Jisugwimundo to some groups. by this way we found the range of magic constant on several HTP.

## 2. Intoduction.

In the hexagonal tortoise problem numbers are arranged into vertices on a graph that is composed of hexagons. Like a normal magic square, a normal HTP uses consecutive integers from 1 to n. The hexagonal sum is the sum of the six numbers on each hexagon, and like magic squares, every hexagonal sum must be same in HTP. Therefore, the magic constant in a HTP is defined as the hexagonal sum. The average magic constant in an n-vertex HTP is $3n+3$.

Theorem 1. If we use an n-vertex m-magic constant HTP, we get an n-vertex $(6n+6-m)$-magic constant HTP.

Proof. Because the magic constant of a HTP is the sum of six numbers, if we give the numbers in reverse order to the HTP, we get a new HTP.

## 3. Diamond-shaped HTP.

A diamond-shaped HTP consists of $n^2$ wedded hexagons. An n by n HTP has $2N^2+4N$ vertices.  If we cut a diamond-shaped HTP in half at n-1 by n-1 and fill it with n linearly-conjoined hexagons, we get a new n by n diamond-shaped HTP. We need 4n+2 more vertices for this transformation. A 1 by 1 HTP consists of 6 vertices. In this way, an n by n HTP has $2N^2+4N$ vertices.

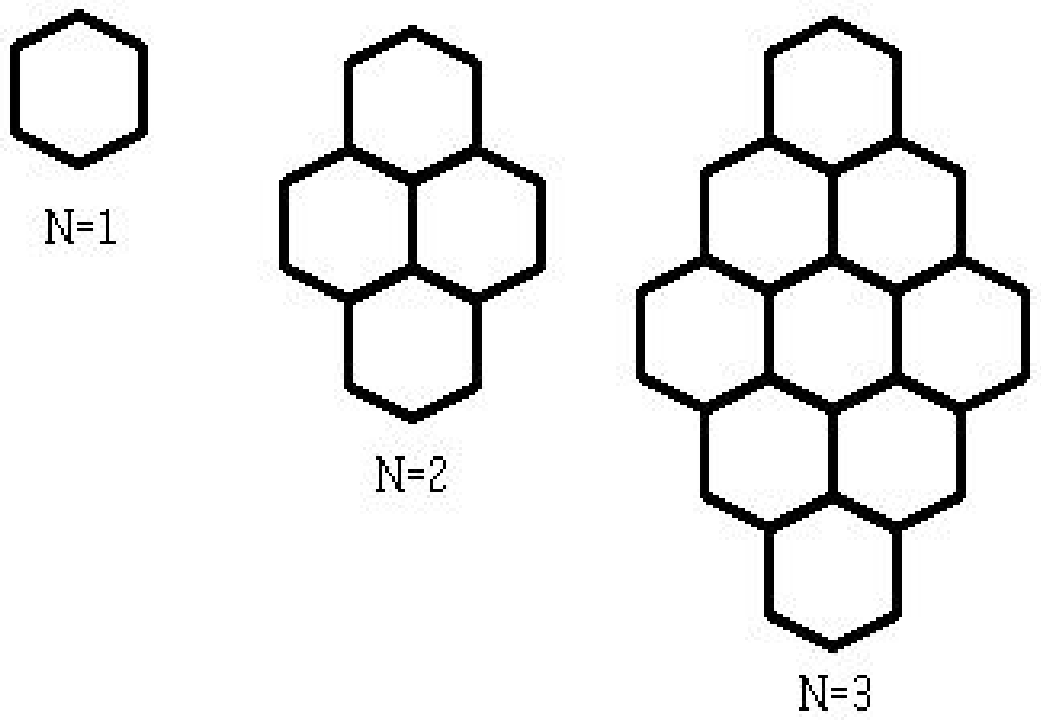

Fig 1 Diamond-shaped HTP

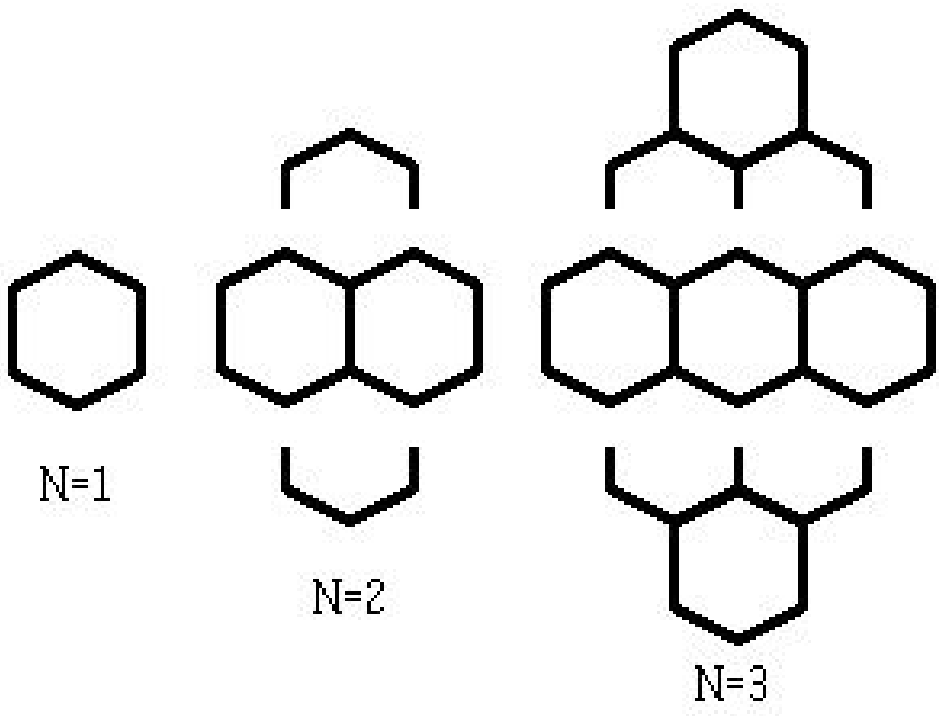

Fig 2 Counting the number of vertices on diamond-shaped HTP

3.1 The range of the magic constant in a diamond-shaped HTP

3.1.0 A 1 by 1 diamond-shaped HTP

A 1 by 1 diamond-shaped HTP consists of six vertices in one hexagon. Its magic constant is obviously 21.

3.1.1 A 2 by 2 diamond-shaped HTP

A 2 by 2 Diamond-shaped HTP has 16 vertices.

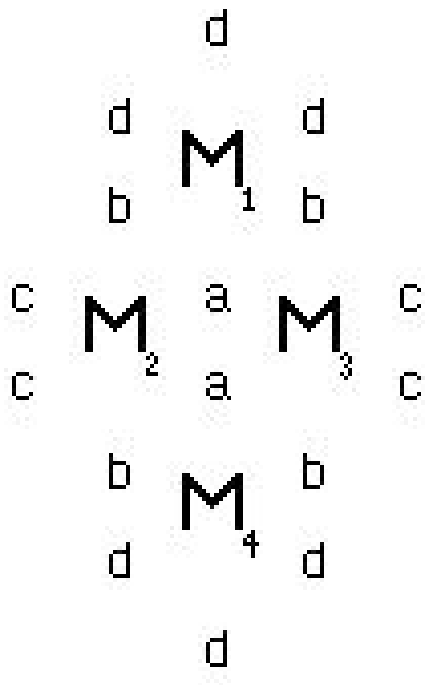

Fig 3

Let us group the 16 vertices into four groups (a, b, c, d) and give name to the four hexagons (M1, M2, M3, M4). This way is depicted in Fig 3. Let the sum of the numbers in the groups a- d be A–D respectively, and the magic constant be M. The sum of the 16 numbers is 136. Then, 136 = A+B+C+D. If we select M1, M4 and c, we choose the 16 numbers once each; therefore, C = 136-2M. On the other hand, If we select M2 and M3, we choose the numbers in group c and b once each, and the numbers in group a twice; 2M = 2A+2B+C. Because C = 136-2M, 4M = 2A+B+136. Therefore, the magic constant can be maximized when A+B and A are individually maximized. A+B is the sum of six numbers in the range 1-16, and A is the sum of two numbers in the range 1-16.

A is maximized when a = {{15,16}}. A+B is maximized when a∪b = {{11,12,13,14,15,16}}. In this case, M = 62. By theorem 1, the minimum magic constant is 40. We will depict whole natural numbers in this range could be a magic constant of 2 by 2 diamond-shaped HTP in Appendix.

3.1.2 A 3 by 3 diamond-shaped HTP

A 3 by 3 diamond-shaped HTP has 30 vertices.

In a 3 by 3 diamond-shaped HTP, it had been shown that the magic constant is in the range 76≤M≤110[2], also 3 by 3 diamond-shaped HTPs that have a magic constant in the range 77≤M≤108 have already been found[2]. By Theorem 1, a 3 by 3 diamond-shaped HTP that has a magic constant of 109 also exists.

We will show that a 3 by 3 diamond-shaped HTP that has a magic constant of 76 or 110 cannot exist. The existence of a 3 by 3 diamond-shaped HTP that has a magic constant of 76 or 110 HTP is equivalent according to theorem 1. Let us show that a 3 by 3 diamond-shaped HTP that has a magic constant of 110 doesn't exist. Let us group 30 vertices into four hexagons and six single vertices G(1), G(2), G(3), G(4), f(1), f(2), f(3), f(4), f(5), f(6). This way is depicted in Fig 4. The hexagonal sums of G(1)–G(4) are M.
4M+f(1)+f(2)+f(3)+f(4)+f(5)+f(6) = 465. Because M = 110, f(1)+f(2)+f(3)+f(4)+f(5)+f(6) = 25. f(2)+f(5) is maximized when the other four numbers f(1), f(3), f(4), and f(6) are in the range 1–4. In this case f(2)+f(5) is 15.

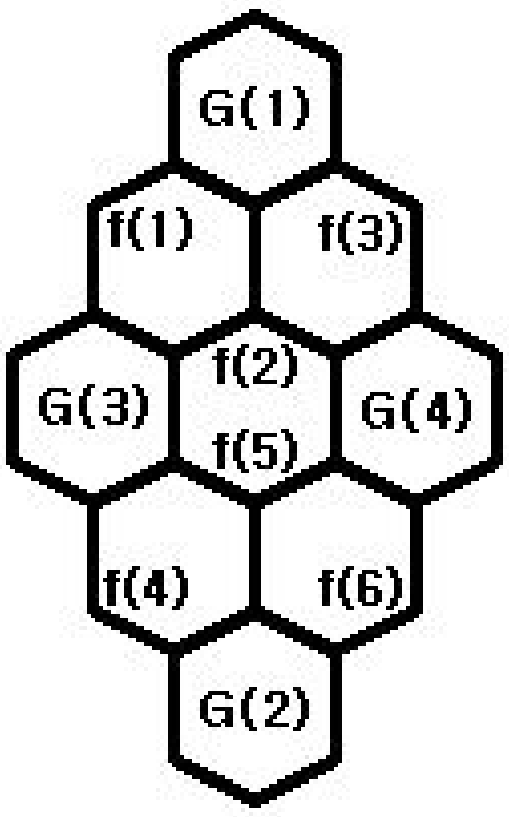

Fig 4 [2]

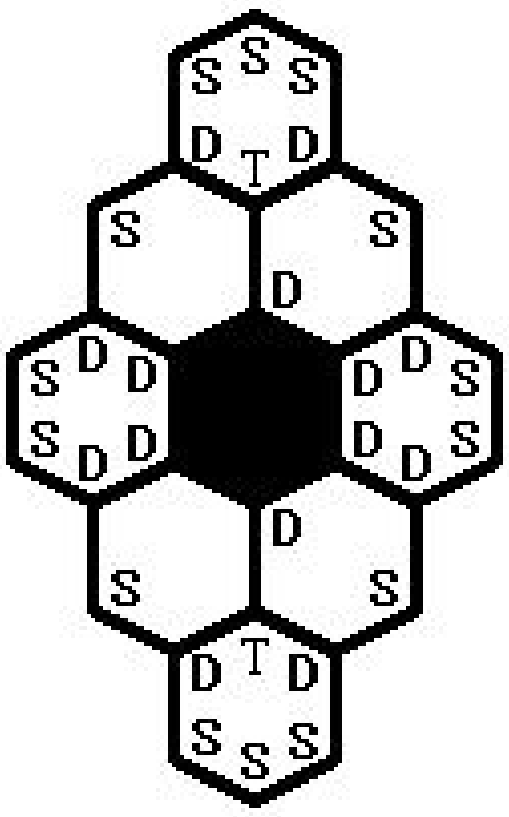

Fig 5 [2]

Let us calculate the sum of the eight hexagonal sums other than the central hexagon; it is 8M. Let us group the 30 vertices by the number of hexagons to which each vertex belongs in the eight hexagons. Two vertices belong to three hexagons, 14 vertices belong to one hexagon, and 14 belong to two hexagons. Let the groups that the vertices belong to three hexagons be T, the vertices that belong to two hexagons be D, and the vertices that belong to one hexagon be S. This way is depicted in Fig 5. Then define the sum of the group T as t, D as d, and S as s. 8M = 3t+2d+s=465+2t+d. 465+2t+d is maximized when t+d and t are individually maximized. However, f(2) and f(5) belong in group D. t+d is maximized when T∪D = {17,18,19,20,21,22,23,24,25,26,27,28,29,30,f(2),f(5)}, and t+d is maximized when T = {29,30}. In this case the maximum of 465+2t+d is 868, which is lower than 880. Therefore, no 3 by 3 diamond-shaped HTP has a magic constant of 110.

We will depict whole natural numbers in this range could be a magic constant of 3 by 3 diamond-shaped HTP in Appendix.

4. Triangular-shaped HTP.

Triangular-shaped HTPs consist of $\frac{n(n+1)}{2}$ wedded hexagons. An n by n HTP has $N^2+4N+1$ vertices. If we add 2n+3 vertices to an n-1 by n-1 diamond-shaped HTP, we get a new n by n triangular-shaped HTP. In this way, an n by n HTP has $N^2+4N+1$ vertices.

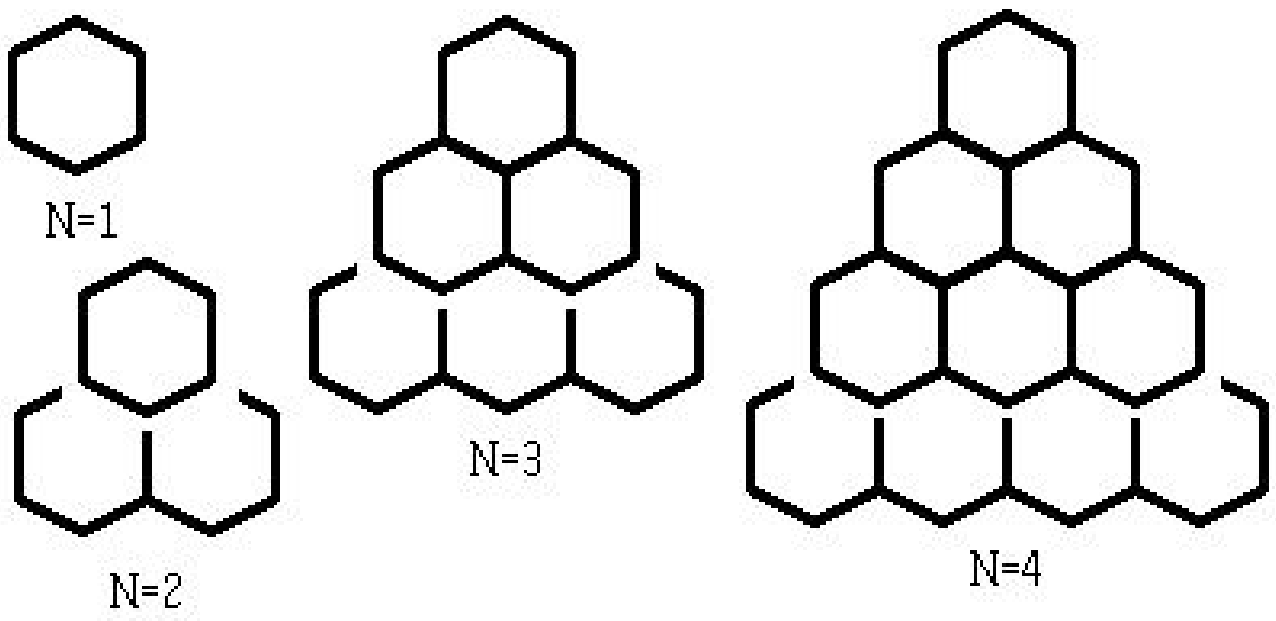

Fig 6 Counting the number of vertices on triangular-shaped HTP

4.1.1 A 2 by 2 triangular-shaped HTP

A 2 by 2 triangular-shaped HTP has 13 vertices. Let us group the 13 vertices in Fig 7.

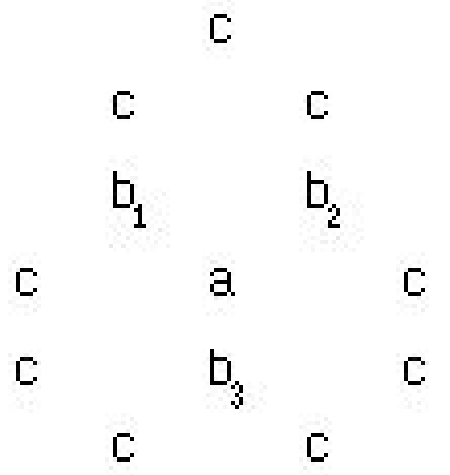

Fig 7

Let us define the sum of the nine numbers in vertex group c as C. The sum of the three hexagonal sums is 3M. 3M = 3a+2(b1+b2+b3)+C = 2a+b1+b2+b3+91. 2a+b1+b2+b3+91 is maximized when a+b1+b2+b3 and a are individually maximized. This can only be possible when a = 13 and {b1,b2,b3} = {10,11,12}. In this case, M = 50. Contrarily, the minimum magic constant is 34. We will depict whole natural numbers in this range could be a magic constant of 2 by 2 triangular-shaped HTP in Appendix.

4.1.2 A 3 by 3 triangular-shaped HTP

A 3 by 3 triangular-shaped HTP has 22 vertices. Let us group the 22 vertices as in Fig 8.

Take all vertices except d1, d2, d3 and a construct three hexagons, and then 253-(a+d1+d2+d3) = 3M. M is maximized when a+d1+d2+d3 is minimized. If {a,d1,d2,d3} = {1,2,3,4}, M is maximized. In this case, M = 71. Meanwhile, the minimum magic constant is 57.

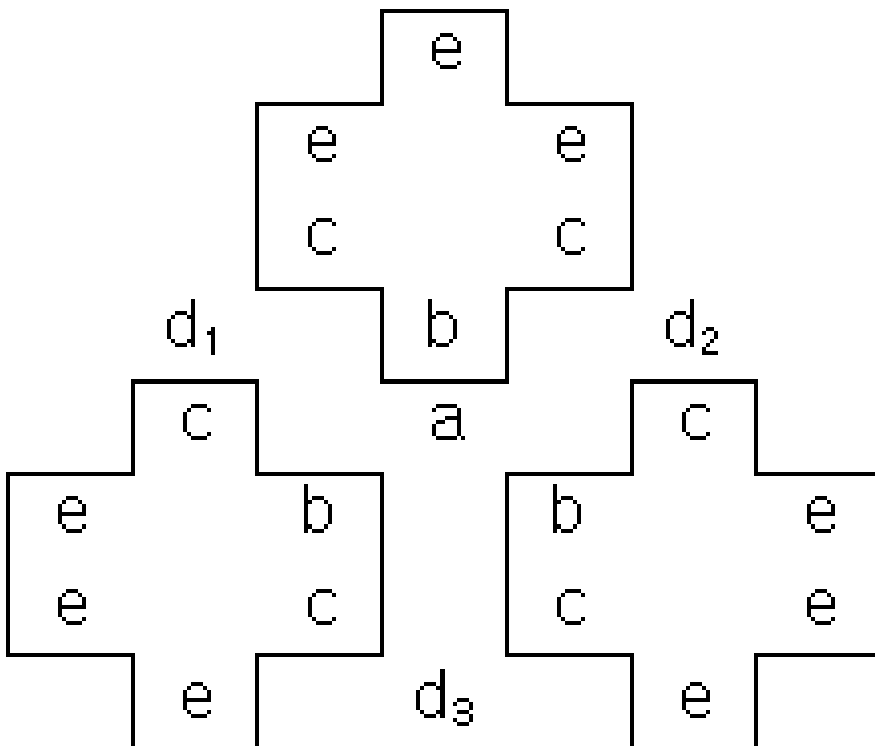

Fig 8 3 by 3 triangular-shaped HTP

We will depict whole natural numbers in this range could be a magic constant of 3 by 3 triangular-shaped HTP in Appendix.

4.1.3 A 4 by 4 triangular-shaped HTP

A 4 by 4 triangular-shaped HTP has 33 vertices. Let give name to the ten hexagons (M1,M2,M3,M4,…,M10), and group the 33 vertices as in Fig 9. Let us define the sum of numbers in vertex group a as A, the sum of b as B, etc. Then a,b,c,d,e,M4 constructs the total set. Therefore, A+B+C+D+E+M = 561. M0,M4,M60,M9,a,b constructs the total set. Therefore, A+B+4M = 561. If we select M1, M2,…, and M10, we choose the numbers in group c and m three times, group d twice, and the numbers in group d once. So 10M = 2A+B+3C+2D+E+3M. Because A+B+C+D+E+M = 561, 8M-561 = 2C+D+A. Because A+B+4M=561, 12M=2C+D-B+1122.

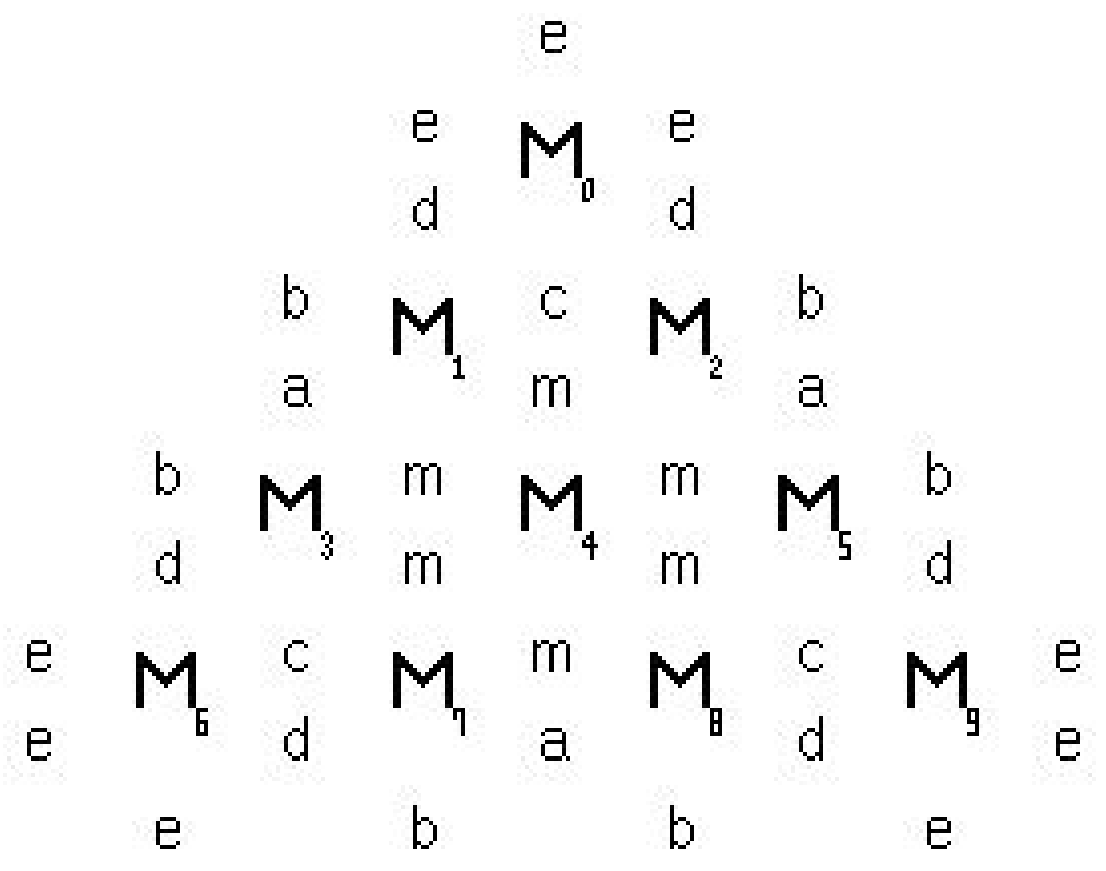

Fig 9 4 by 4 triangular-shaped HTP

Then 2C+D-B+1122 is maximized when C and C+D are individually maximized and b is individually minimized. This can only be possible when c = {33,32,31}, d = {30,29,28,27,26,25}, and b = {1,2,3,4,5,6}. In this case M = 121.5, so M cannot be bigger than 122. In addition, the lower bound of M is 83. We will depict whole natural numbers in this range could be a magic constant of 4 by 4 triangular-shaped HTP in Appendix.

5. Hexagonal-shaped HTP.

A hexagonal-shaped HTP consists of 3n(n-1)+1 wedded hexagons (Centered hexagonal number). An n by n HTP has $6N^2$ vertices. If we add one layer to an n-1 by n-1 hexagonal-shaped HTP, we get a new n by n hexagonal-shaped HTP. We need 12n-6 more vertices to perform this transformation. A 1 by 1 HTP consists of 6 vertices. In this way, an n by n HTP has $6N^2$ vertices.

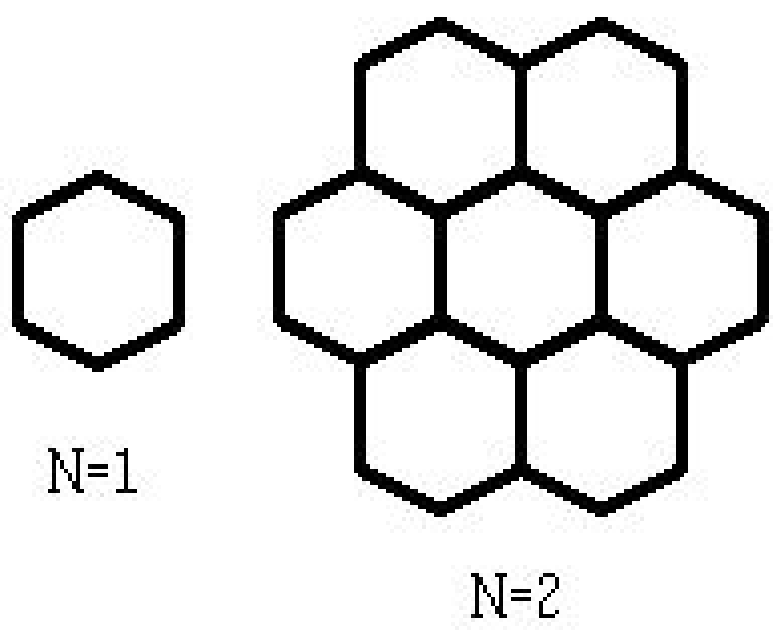

Fig 10 hexagonal-shaped HTP

5.1 The range of the magic constant in hexagonal-shaped HTPs

5.1.0 A 1 by 1 hexagonal-shaped HTP

A 1 by 1 hexagonal-shaped HTP consists of six vertices in one hexagon. Its magic constant is obviously 21.

### 5.1.1 A 2 by 2 hexagonal-shaped HTP

The 2 by 2 hexagonal-shaped HTP has 24 vertices. Let us group these 24 vertices as in Fig 11 and give names to the seven hexagons (M0,M1,M2,M3,M4,M5,M6). Let calculate the sum of the six hexagonal sums of M1–M6, it is 6M. If we select M1 M2,…, and M6, we choose all the numbers in group b and m two times, and those in group c once. So 6M = 2M+2B+C. Because M+B+C = 300, 300+B = 5M. Then 300+B is maximized when B is maximized. It can be possible only when b = {19,20,21,22,23,24}. In this case M=85.8, So M cannot be bigger than 85. In addition, the lower bound of M is 65.

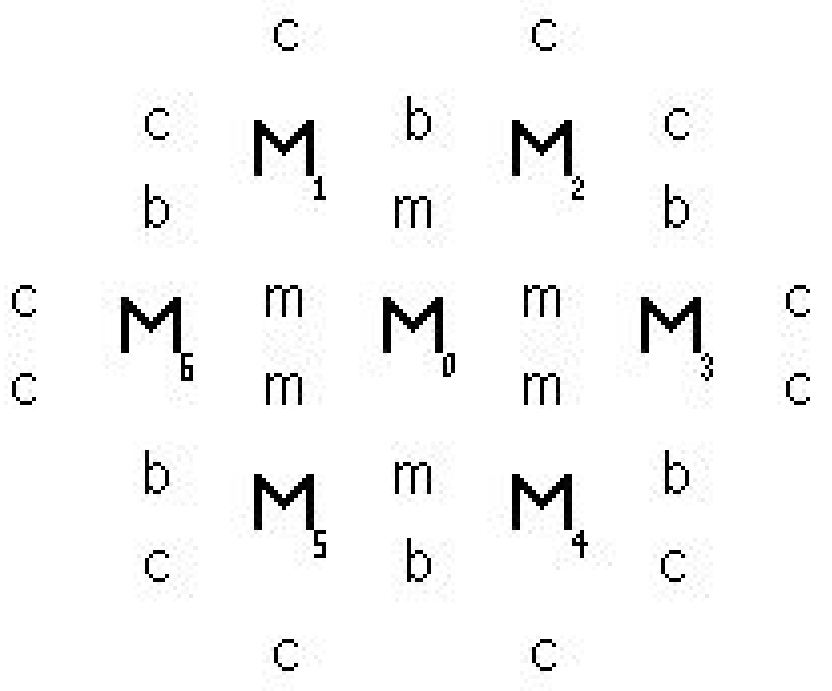

Fig 11

### 5.1.2 A 3 by 3 hexagonal-shaped HTP

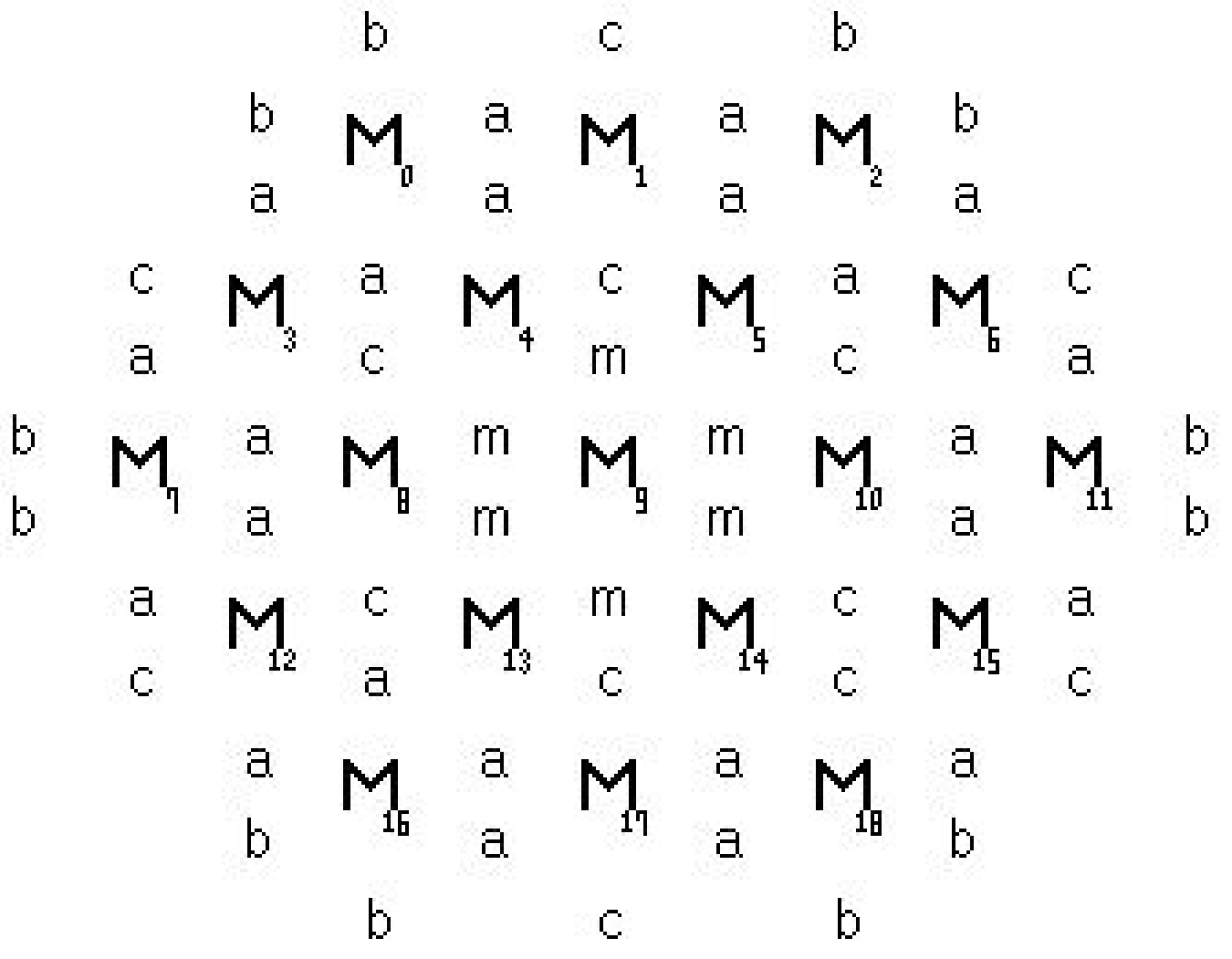

Fig 12 3 by 3 hexagonal-shaped HTP

A 3 by 3 hexagonal-shaped HTP has 54 vertices. Let us give names to ten hexagons (M0,M1,M2,M3,M4,…,M18) and group the 54 vertices as in Fig 12. Let us define the sum of numbers in the vertex group a as A, the sum of b as B, etc., then a,c,m constructs seven hexagons. In addition , a,b,m constructs seven hexagons. So A+C+M = A+B+M = 7M. Because A+B+C+M = 1485, B = C = 1485-7M. B+C = 2970-14M. M is maximized when B+C is minimized. In this case, b∪c = {1,2,…,24} and B+C = 300. The Upper bound of M = 190. In the same manner, the lower bound of M = 140. We will depict whole natural numbers in this range could be a magic constant of 3 by 3 hexagonal-shaped HTP in Appendix.

## 6. A 2 by 2 hexagonal star-shaped HTP.

A 2 by 2 hexagonal star-shaped HTP has 42 vertices. We can pick seven hexagons from the HTP and choose all vertices at once. Therefore, the magic constant must be 129.

Appendix

3.1.1 2 by 2 diamond-shaped HTP

| | | | | | | | | | | | | | | | | | | | |
|---|---|---|---|---|---|---|---|---|---|---|---|---|---|---|---|---|---|---|---|
| 62 | | 8 | | | 61 | | 5 | | | 60 | | 9 | | | 59 | | 5 | | |
| | 3 | | 9 | | | 4 | | 7 | | | 6 | | 12 | | | 4 | | 15 | |
| | 12 | | 14 | | | 14 | | 15 | | | 7 | | 10 | | | 9 | | 10 | |
| 2 | | 16 | | 1 | 1 | | 16 | | 2 | 1 | | 16 | | 2 | 1 | | 16 | | 2 |
| 4 | | 15 | | 5 | 8 | | 13 | | 3 | 8 | | 13 | | 5 | 8 | | 13 | | 7 |
| | 13 | | 11 | | | 9 | | 12 | | | 15 | | 14 | | | 12 | | 11 | |
| | 10 | | 6 | | | 6 | | 11 | | | 3 | | 11 | | | 3 | | 14 | |
| | | 7 | | | | | 10 | | | | | 4 | | | | | 6 | | |
| 58 | | 13 | | | 57 | | 3 | | | 56 | | 14 | | | 55 | | 11 | | |
| | 1 | | 10 | | | 12 | | 15 | | | 3 | | 13 | | | 1 | | 12 | |
| | 8 | | 11 | | | 11 | | 2 | | | 10 | | 1 | | | 8 | | 9 | |
| 3 | | 15 | | 4 | 4 | | 14 | | 5 | 7 | | 15 | | 5 | 15 | | 14 | | 5 |
| 6 | | 12 | | 7 | 6 | | 13 | | 7 | 4 | | 11 | | 8 | 2 | | 13 | | 4 |
| | 14 | | 9 | | | 9 | | 16 | | | 9 | | 16 | | | 3 | | 10 | |
| | 16 | | 5 | | | 1 | | 10 | | | 2 | | 12 | | | 6 | | 7 | |
| | | 2 | | | | | 8 | | | | | 6 | | | | | 16 | | |
| 54 | | 8 | | | 53 | | 8 | | | 52 | | 13 | | | 51 | | 3 | | |
| | 5 | | 16 | | | 6 | | 16 | | | 2 | | 11 | | | 10 | | 11 | |
| | 1 | | 10 | | | 1 | | 9 | | | 3 | | 8 | | | 5 | | 6 | |
| 15 | | 14 | | 7 | 15 | | 13 | | 10 | 12 | | 15 | | 9 | 13 | | 16 | | 4 |
| 2 | | 13 | | 4 | 3 | | 14 | | 2 | 5 | | 10 | | 6 | 2 | | 1 | | 15 |
| | 9 | | 6 | | | 7 | | 5 | | | 7 | | 4 | | | 14 | | 9 | |
| | 11 | | 12 | | | 4 | | 12 | | | 1 | | 16 | | | 7 | | 8 | |
| | | 3 | | | | | 11 | | | | | 14 | | | | | 12 | | |
| 50 | | 4 | | | 49 | | 9 | | | 48 | | 9 | | | 47 | | 6 | | |
| | 15 | | 6 | | | 11 | | 1 | | | 12 | | 1 | | | 16 | | 5 | |
| | 14 | | 9 | | | 16 | | 8 | | | 16 | | 7 | | | 9 | | 8 | |
| 5 | | 2 | | 8 | 2 | | 4 | | 7 | 2 | | 3 | | 10 | 2 | | 3 | | 12 |
| 12 | | 7 | | 11 | 14 | | 3 | | 15 | 15 | | 4 | | 13 | 15 | | 4 | | 13 |
| | 10 | | 13 | | | 10 | | 12 | | | 8 | | 11 | | | 14 | | 7 | |
| | 16 | | 1 | | | 13 | | 5 | | | 6 | | 5 | | | 11 | | 10 | |
| | | 3 | | | | | 6 | | | | | 14 | | | | | 1 | | |
| 46 | | 3 | | | 45 | | 14 | | | 44 | | 4 | | | 43 | | 12 | | |
| | 14 | | 4 | | | 5 | | 2 | | | 16 | | 7 | | | 13 | | 2 | |
| | 7 | | 16 | | | 6 | | 15 | | | 9 | | 6 | | | 8 | | 7 | |
| 10 | | 2 | | 12 | 13 | | 3 | | 12 | 14 | | 2 | | 13 | 16 | | 1 | | 15 |
| 13 | | 6 | | 9 | 11 | | 4 | | 10 | 11 | | 5 | | 10 | 9 | | 4 | | 10 |
| | 8 | | 1 | | | 8 | | 1 | | | 3 | | 8 | | | 5 | | 6 | |
| | 15 | | 5 | | | 16 | | 7 | | | 1 | | 12 | | | 14 | | 3 | |
| | | 11 | | | | | 9 | | | | | 15 | | | | | 11 | | |

| | | | | | | | | | | | | | | |
|---|---|---|---|---|---|---|---|---|---|---|---|---|---|---|
| 42 | | 8 | | | 41 | | 12 | | | 40 | | 9 | | |
| | 11 | | 5 | | | 13 | | 10 | | | 14 | | 8 | |
| | 10 | | 7 | | | 3 | | 2 | | | 5 | | 3 | |
| 16 | | 1 | | 15 | 16 | | 1 | | 15 | 15 | | 1 | | 16 |
| 9 | | 4 | | 12 | 9 | | 4 | | 14 | 13 | | 2 | | 12 |
| | 2 | | 3 | | | 8 | | 5 | | | 4 | | 6 | |
| | 14 | | 6 | | | 11 | | 6 | | | 7 | | 11 | |
| | | 13 | | | | | 7 | | | | | 10 | | |

3.1.2 3 by 3 diamond-shaped HTP

| | | | | | | | | | | | | | | | | | | | | |
|---|---|---|---|---|---|---|---|---|---|---|---|---|---|---|---|---|---|---|---|---|
| | | | 19 | | | | 78 | | | 23 | | | | 79 | | | 25 | | | |
| 77 | | 24 | | 25 | | | | | 24 | | 25 | | | | | 24 | | 23 | | |
| | | 4 | | 3 | | | | | 3 | | 2 | | | | | 1 | | 2 | | |
| | 28 | | 2 | | 27 | | | 27 | | 1 | | 30 | | | 30 | | 4 | | 27 | |
| | 6 | | 26 | | 9 | | | 16 | | 26 | | 4 | | | 6 | | 26 | | 13 | |
| 20 | | 11 | | 10 | | 22 | 21 | | 5 | | 15 | | 19 | 16 | | 12 | | 7 | | 14 |
| 23 | | 5 | | 8 | | 21 | 22 | | 8 | | 11 | | 20 | 17 | | 10 | | 15 | | 11 |
| | 12 | | 17 | | 7 | | | 6 | | 13 | | 9 | | | 18 | | 9 | | 19 | |
| | 29 | | 1 | | 30 | | | 29 | | 10 | | 28 | | | 29 | | 5 | | 28 | |
| | | 13 | | 14 | | | | | 12 | | 7 | | | | | 8 | | 3 | | |
| | | 18 | | 16 | | | | | 18 | | 17 | | | | | 22 | | 20 | | |
| | | | 15 | | | | | | | 14 | | | | | | | 21 | | | |
| 80 | | | 17 | | | | 81 | | | 24 | | | | 82 | | | 16 | | | |
| | | 18 | | 14 | | | | | 9 | | 7 | | | | | 7 | | 13 | | |
| | | 11 | | 12 | | | | | 18 | | 22 | | | | | 24 | | 14 | | |
| | 28 | | 8 | | 26 | | | 29 | | 1 | | 25 | | | 26 | | 8 | | 28 | |
| | 6 | | 5 | | 9 | | | 2 | | 16 | | 6 | | | 1 | | 4 | | 17 | |
| 24 | | 22 | | 20 | | 21 | 21 | | 15 | | 11 | | 27 | 23 | | 19 | | 11 | | 21 |
| 25 | | 2 | | 4 | | 23 | 26 | | 14 | | 12 | | 20 | 27 | | 10 | | 18 | | 12 |
| | 1 | | 27 | | 3 | | | 3 | | 13 | | 5 | | | 2 | | 20 | | 3 | |
| | 30 | | 7 | | 29 | | | 28 | | 4 | | 30 | | | 29 | | 6 | | 30 | |
| | | 13 | | 10 | | | | | 19 | | 17 | | | | | 15 | | 5 | | |
| | | 19 | | 16 | | | | | 10 | | 8 | | | | | 25 | | 9 | | |
| | | | 15 | | | | | | | 23 | | | | | | | 22 | | | |
| 83 | | | 17 | | | | 84 | | | 4 | | | | 85 | | | 25 | | | |
| | | 25 | | 6 | | | | | 27 | | 11 | | | | | 16 | | 9 | | |
| | | 22 | | 12 | | | | | 22 | | 19 | | | | | 21 | | 12 | | |
| | 24 | | 1 | | 29 | | | 25 | | 1 | | 26 | | | 30 | | 2 | | 29 | |
| | 9 | | 8 | | 13 | | | 8 | | 5 | | 13 | | | 4 | | 5 | | 11 | |
| 16 | | 19 | | 20 | | 23 | 24 | | 23 | | 20 | | 7 | 22 | | 23 | | 26 | | 3 |
| 26 | | 11 | | 10 | | 14 | 15 | | 12 | | 10 | | 28 | 14 | | 7 | | 18 | | 8 |
| | 2 | | 15 | | 3 | | | 2 | | 14 | | 6 | | | 15 | | 6 | | 19 | |
| | 30 | | 7 | | 27 | | | 30 | | 9 | | 29 | | | 27 | | 13 | | 28 | |
| | | 18 | | 21 | | | | | 17 | | 16 | | | | | 17 | | 1 | | |
| | | 4 | | 28 | | | | | 18 | | 3 | | | | | 10 | | 24 | | |
| | | | 5 | | | | | | | 21 | | | | | | | 20 | | | |
| 86 | | | 19 | | | | 87 | | | 14 | | | | 88 | | | 27 | | | |
| | | 2 | | 22 | | | | | 24 | | 8 | | | | | 21 | | 5 | | |
| | | 18 | | 20 | | | | | 11 | | 23 | | | | | 3 | | 4 | | |
| | 28 | | 5 | | 27 | | | 28 | | 7 | | 27 | | | 18 | | 28 | | 24 | |
| | 9 | | 1 | | 7 | | | 15 | | 1 | | 3 | | | 17 | | 10 | | 9 | |
| 10 | | 25 | | 26 | | 16 | 5 | | 25 | | 26 | | 19 | 26 | | 12 | | 13 | | 14 |
| 17 | | 13 | | 15 | | 14 | 6 | | 16 | | 17 | | 4 | 2 | | 23 | | 11 | | 25 |
| | 12 | | 6 | | 8 | | | 20 | | 2 | | 18 | | | 8 | | 19 | | 16 | |
| | 30 | | 4 | | 29 | | | 30 | | 9 | | 29 | | | 22 | | 15 | | 20 | |
| | | 21 | | 24 | | | | | 10 | | 12 | | | | | 1 | | 7 | | |
| | | 23 | | 3 | | | | | 21 | | 13 | | | | | 30 | | 29 | | |
| | | | 11 | | | | | | | 22 | | | | | | | 6 | | | |
| 89 | | | 1 | | | | 90 | | | 29 | | | | 91 | | | 1 | | | |
| | | 30 | | 24 | | | | | 1 | | 3 | | | | | 30 | | 27 | | |
| | | 3 | | 9 | | | | | 28 | | 27 | | | | | 2 | | 5 | | |
| | 28 | | 22 | | 16 | | | 14 | | 2 | | 15 | | | 29 | | 26 | | 23 | |
| | 5 | | 11 | | 17 | | | 4 | | 17 | | 8 | | | 3 | | 6 | | 9 | |

| | | | | | | |
|---|---|---|---|---|---|---|
| 26 | | 20 | | 14 | | 8 |
| 7 | | 13 | | 19 | | 25 |
| | 18 | | 12 | | 6 | |
| | 15 | | 21 | | 27 | |
| | | 10 | | 4 | | |
| | | 23 | | 29 | | |
| | | | 2 | | | |
| 26 | | 25 | | 21 | | 22 |
| 18 | | 5 | | 9 | | 7 |
| | 12 | | 13 | | 23 | |
| | 30 | | 10 | | 16 | |
| | | 20 | | 19 | | |
| | | 24 | | 11 | | |
| | | | 6 | | | |
| 28 | | 25 | | 22 | | 19 |
| 4 | | 7 | | 10 | | 13 |
| | 24 | | 21 | | 18 | |
| | 8 | | 11 | | 14 | |
| | | 20 | | 17 | | |
| | | 12 | | 15 | | |
| | | | 16 | | | |
| 92 | | | 30 | | | |
| | | 22 | | 4 | | |
| | | 24 | | 9 | | |
| | 10 | | 3 | | 18 | |
| | 21 | | 27 | | 19 | |
| 26 | | 7 | | 16 | | 15 |
| 12 | | 20 | | 17 | | 11 |
| | 6 | | 5 | | 14 | |
| | 8 | | 25 | | 29 | |
| | | 28 | | 2 | | |
| | | 23 | | 13 | | |
| | | | 1 | | | |
| 93 | | | 5 | | | |
| | | 26 | | 4 | | |
| | | 10 | | 27 | | |
| | 24 | | 21 | | 9 | |
| | 6 | | 15 | | 8 | |
| 3 | | 17 | | 13 | | 20 |
| 28 | | 14 | | 18 | | 11 |
| | 25 | | 16 | | 23 | |
| | 7 | | 12 | | 22 | |
| | | 19 | | 2 | | |
| | | 30 | | 29 | | |
| | | | 1 | | | |

### 5.1.1 2 by 2 hexagonal-shaped HTP

| | | | | | | |
|---|---|---|---|---|---|---|
| 65 | | 19 | | 11 | | |
| | 16 | | 3 | | 24 | |
| | 10 | | 8 | | 2 | |
| 20 | | 9 | | 17 | | 12 |
| 15 | | 6 | | 7 | | 23 |
| | 5 | | 18 | | 4 | |
| | 21 | | 1 | | 13 | |
| | | 14 | | 22 | | |
| 66 | | 20 | | 16 | | |
| | 23 | | 2 | | 18 | |
| | 5 | | 9 | | 4 | |
| 10 | | 7 | | 17 | | 13 |
| 24 | | 19 | | 8 | | 21 |
| | 1 | | 6 | | 3 | |
| | 11 | | 15 | | 12 | |
| | | 14 | | 22 | | |
| 67 | | 20 | | 16 | | |
| | 24 | | 1 | | 14 | |
| | 10 | | 7 | | 8 | |
| 19 | | 5 | | 21 | | 9 |
| 13 | | 18 | | 12 | | 6 |
| | 2 | | 4 | | 11 | |
| | 17 | | 3 | | 22 | |
| | | 23 | | 15 | | |
| 68 | | 20 | | 23 | | |
| | 13 | | 7 | | 11 | |
| | 3 | | 10 | | 1 | |
| 5 | | 15 | | 16 | | 24 |
| 21 | | 6 | | 17 | | 8 |
| | 18 | | 4 | | 2 | |
| | 19 | | 9 | | 22 | |
| | | 12 | | 14 | | |
| 69 | | 9 | | 24 | | |
| | 21 | | 3 | | 12 | |
| | 11 | | 10 | | 4 | |
| 1 | | 15 | | 16 | | 22 |
| 19 | | 5 | | 17 | | 8 |
| | 18 | | 6 | | 2 | |
| | 20 | | 7 | | 14 | |
| | | 13 | | 23 | | |
| 70 | | 10 | | 1 | | |
| | 21 | | 13 | | 22 | |
| | 5 | | 18 | | 7 | |
| 23 | | 3 | | 9 | | 14 |
| 16 | | 8 | | 12 | | 24 |
| | 15 | | 20 | | 4 | |
| | 2 | | 6 | | 11 | |
| | | 19 | | 17 | | |
| 71 | | 10 | | 23 | | |
| | 22 | | 13 | | 1 | |
| | 5 | | 18 | | 7 | |
| 24 | | 3 | | 9 | | 20 |
| 16 | | 8 | | 12 | | 19 |
| | 15 | | 21 | | 4 | |
| | 2 | | 11 | | 6 | |
| | | 14 | | 17 | | |
| 72 | | 24 | | 20 | | |
| | 6 | | 15 | | 2 | |
| | 5 | | 19 | | 7 | |
| 16 | | 3 | | 9 | | 23 |
| 22 | | 8 | | 12 | | 17 |
| | 18 | | 21 | | 4 | |
| | 13 | | 11 | | 10 | |
| | | 1 | | 14 | | |
| 73 | | 17 | | 21 | | |
| | 22 | | 7 | | 9 | |
| | 6 | | 1 | | 12 | |
| 19 | | 20 | | 23 | | 8 |
| 10 | | 2 | | 3 | | 14 |
| | 16 | | 24 | | 13 | |
| | 5 | | 11 | | 4 | |
| | | 15 | | 18 | | |
| 74 | | 14 | | 2 | | |
| | 12 | | 13 | | 24 | |
| | 9 | | 11 | | 1 | |
| 16 | | 15 | | 23 | | 4 |
| 10 | | 5 | | 3 | | 22 |
| | 19 | | 17 | | 21 | |
| | 8 | | 7 | | 6 | |
| | | 18 | | 20 | | |
| 75 | | 15 | | 14 | | |
| | 10 | | 22 | | | |
| | 4 | | 3 | | 11 | |
| 13 | | 21 | | 5 | 20 | |
| 12 | | 6 | | 23 | | 24 |
| | 19 | | 17 | | | 1 |
| | 18 | | 8 | | 2 | |
| | | 7 | | 9 | 16 | |

5.1.2 3 by 3 hexagonal-shaped HTP

| | | | | | | | | | | |
|---|---|---|---|---|---|---|---|---|---|---|
| 190 | | | 16 | | 11 | | 22 | | | |
| | | 13 | | 24 | | 35 | | 6 | | |
| | | 50 | | 53 | | 48 | | 33 | | |
| | 2 | | 34 | | 19 | | 46 | | 15 | |
| | 31 | | 32 | | 23 | | 12 | | 39 | |
| 14 | | 41 | | 29 | | 42 | | 45 | | 21 |
| 3 | | 52 | | 28 | | 30 | | 54 | | 4 |
| | 49 | | 8 | | 38 | | 7 | | 27 | |
| | 5 | | 40 | | 25 | | 47 | | 18 | |
| | | 36 | | 51 | | 43 | | 37 | | |
| | | 10 | | 44 | | 26 | | 17 | | |
| | | | 9 | | 1 | | 20 | | | |
| 189 | | | 17 | | 13 | | 24 | | | |
| | | 7 | | 50 | | 30 | | 12 | | |
| | | 29 | | 47 | | 48 | | 37 | | |
| | 9 | | 39 | | 1 | | 38 | | 3 | |
| | 34 | | 27 | | 53 | | 28 | | 31 | |
| 16 | | 51 | | 22 | | 21 | | 52 | | 23 |
| 11 | | 45 | | 19 | | 20 | | 42 | | 6 |
| | 32 | | 25 | | 54 | | 26 | | 35 | |
| | 14 | | 40 | | 2 | | 46 | | 4 | |
| | | 33 | | 49 | | 41 | | 36 | | |
| | | 18 | | 44 | | 43 | | 8 | | |
| | | | 5 | | 10 | | 15 | | | |
| 188 | | | 12 | | 6 | | 30 | | | |
| | | 16 | | 46 | | 24 | | 17 | | |
| | | 20 | | 41 | | 40 | | 23 | | |
| | 10 | | 53 | | 31 | | 54 | | 15 | |
| | 50 | | 3 | | 26 | | 1 | | 44 | |
| 13 | | 52 | | 34 | | 36 | | 51 | | 18 |
| 14 | | 38 | | 28 | | 29 | | 39 | | 11 |
| | 21 | | 33 | | 35 | | 32 | | 25 | |
| | 5 | | 42 | | 2 | | 43 | | 22 | |
| | | 49 | | 48 | | 47 | | 27 | | |
| | | 4 | | 37 | | 45 | | 29 | | |
| | | | 8 | | 9 | | 7 | | | |
| 187 | | | 7 | | 4 | | 8 | | | |
| | | 17 | | 30 | | 31 | | 14 | | |
| | | 32 | | 50 | | 44 | | 38 | | |
| | 22 | | 51 | | 28 | | 52 | | 12 | |
| | 33 | | 10 | | 2 | | 18 | | 42 | |
| 24 | | 39 | | 46 | | 43 | | 25 | | 23 |
| 20 | | 36 | | 47 | | 48 | | 40 | | 16 |
| | 35 | | 9 | | 1 | | 13 | | 41 | |
| | 19 | | 54 | | 27 | | 53 | | 3 | |
| | | 34 | | 49 | | 45 | | 37 | | |
| | | 6 | | 29 | | 26 | | 21 | | |
| | | | 15 | | 11 | | 5 | | | |
| 186 | | | 16 | | 1 | | 4 | | | |
| | | 15 | | 31 | | 46 | | 20 | | |
| | | 35 | | 50 | | 51 | | 27 | | |
| | 3 | | 39 | | 7 | | 38 | | 18 | |
| | 44 | | 28 | | 8 | | 29 | | 40 | |
| 5 | | 37 | | 54 | | 53 | | 34 | | 23 |
| 21 | | 47 | | 9 | | 10 | | 48 | | 19 |
| | 32 | | 11 | | 52 | | 12 | | 22 | |
| | 25 | | 45 | | 33 | | 49 | | 14 | |
| | | 26 | | 36 | | 30 | | 41 | | |
| | | 24 | | 42 | | 43 | | 6 | | |
| | | | 13 | | 2 | | 17 | | | |
| 185 | | | 11 | | 21 | | 20 | | | |
| | | 13 | | 27 | | 7 | | 16 | | |
| | | 31 | | 51 | | 50 | | 38 | | |
| | 18 | | 52 | | 29 | | 54 | | 14 | |
| | 34 | | 4 | | 6 | | 2 | | 32 | |
| 1 | | 46 | | 43 | | 44 | | 45 | | 25 |
| 17 | | 48 | | 41 | | 42 | | 40 | | 10 |
| | 39 | | 3 | | 9 | | 22 | | 33 | |
| | 5 | | 53 | | 30 | | 35 | | 19 | |
| | | 37 | | 49 | | 47 | | 36 | | |
| | | 26 | | 8 | | 28 | | 24 | | |
| | | | 12 | | 23 | | 15 | | | |
| 184 | | | 6 | | 22 | | 15 | | | |
| | | 11 | | 29 | | 2 | | 21 | | |
| | | 34 | | 51 | | 52 | | 40 | | |
| | 14 | | 53 | | 28 | | 54 | | 12 | |
| | 37 | | 8 | | 1 | | 3 | | 39 | |
| 23 | | 38 | | 43 | | 46 | | 36 | | 16 |
| 10 | | 41 | | 45 | | 44 | | 42 | | 20 |
| | 35 | | 9 | | 5 | | 13 | | 31 | |
| | 18 | | 48 | | 27 | | 49 | | 17 | |
| | | 33 | | 50 | | 47 | | 32 | | |
| | | 24 | | 4 | | 30 | | 19 | | |
| | | | 25 | | 26 | | 7 | | | |
| 183 | | | 26 | | 22 | | 28 | | | |
| | | 8 | | 13 | | 42 | | 11 | | |
| | | 48 | | 36 | | 37 | | 14 | | |
| | 7 | | 52 | | 33 | | 51 | | 25 | |
| | 21 | | 1 | | 16 | | 2 | | 41 | |
| 24 | | 54 | | 45 | | 44 | | 50 | | 29 |
| 10 | | 34 | | 17 | | 18 | | 38 | | 6 |
| | 40 | | 32 | | 43 | | 31 | | 19 | |
| | 27 | | 35 | | 3 | | 39 | | 9 | |
| | | 15 | | 53 | | 49 | | 47 | | |
| | | 30 | | 46 | | 20 | | 23 | | |
| | | | 4 | | 12 | | 5 | | | |
| 182 | | | 12 | | 27 | | 30 | | | |
| | | 25 | | 36 | | 5 | | 20 | | |
| | | 21 | | 41 | | 40 | | 39 | | |
| | 9 | | 47 | | 33 | | 48 | | 24 | |
| | 51 | | 1 | | 15 | | 2 | | 19 | |
| 6 | | 53 | | 45 | | 44 | | 50 | | 22 |
| 28 | | 34 | | 17 | | 18 | | 37 | | 8 |
| | 10 | | 32 | | 43 | | 31 | | 46 | |
| | 29 | | 35 | | 3 | | 38 | | 7 | |
| | | 42 | | 52 | | 49 | | 23 | | |
| | | 16 | | 11 | | 54 | | 14 | | |
| | | | 26 | | 13 | | 4 | | | |
| 181 | | | 12 | | 21 | | 32 | | | |
| | | 13 | | 1 | | 36 | | 11 | | |
| | | 52 | | 49 | | 45 | | 7 | | |
| | 22 | | 54 | | 29 | | 50 | | 24 | |
| | 2 | | 3 | | 5 | | 9 | | 38 | |
| 27 | | 48 | | 41 | 181 | 43 | | 53 | | 19 |
| 28 | | 37 | | 42 | | 6 | | 40 | | 8 |
| | 39 | | 10 | | 44 | | 30 | | 23 | |
| | 26 | | 34 | | 4 | | 46 | | 25 | |
| | | 35 | | 47 | | 51 | | 17 | | |
| | | 14 | | 33 | | 31 | | 16 | | |
| | | | 18 | | 15 | | 20 | | | |
| 180 | | | 13 | | 35 | | 11 | | | |
| | | 24 | | 28 | | 21 | | 18 | | |
| | | 27 | | 51 | | 42 | | 48 | | |
| | 6 | | 37 | | 3 | | 40 | | 7 | |
| | 29 | | 31 | | 54 | | 32 | | 17 | |
| 23 | | 50 | | 4 | | 9 | | 36 | | 14 |
| 25 | | 41 | | 52 | | 53 | | 49 | | 19 |
| | 12 | | 2 | | 8 | | 1 | | 45 | |
| | 44 | | 47 | | 33 | | 39 | | 26 | |
| | | 34 | | 38 | | 46 | | 20 | | |
| | | 16 | | 15 | | 43 | | 10 | | |
| | | | 30 | | 5 | | 22 | | | |
| 179 | | | 14 | | 32 | | 30 | | | |
| | | 15 | | 2 | | 47 | | 25 | | |
| | | 44 | | 54 | | 39 | | 1 | | |
| | 19 | | 50 | | 5 | | 37 | | 24 | |
| | 6 | | 9 | | 35 | | 36 | | 43 | |
| 18 | | 51 | | 26 | | 27 | | 38 | | 31 |
| 12 | | 52 | | 28 | | 29 | | 41 | | 23 |
| | 40 | | 13 | | 34 | | 8 | | 3 | |
| | 21 | | 46 | | 10 | | 45 | | 33 | |
| | | 7 | | 48 | | 53 | | 49 | | |
| | | 16 | | 42 | | 4 | | 11 | | |
| | | | 20 | | 22 | | 17 | | | |
| 178 | | | 14 | | 10 | | 15 | | | |
| | | 23 | | 38 | | 5 | | 16 | | |
| | | 2 | | 47 | | 45 | | 46 | | |
| | 26 | | 54 | | 33 | | 51 | | 25 | |
| | 39 | | 9 | | 1 | | 13 | | 3 | |
| 17 | | 48 | | 34 | | 35 | | 40 | | 20 |
| 19 | | 49 | | 27 | | 28 | | 50 | | 21 |
| | 6 | | 11 | | 53 | | 12 | | 44 | |
| | 32 | | 37 | | 8 | | 36 | | 29 | |
| | | 43 | | 42 | | 41 | | 7 | | |
| | | 22 | | 4 | | 52 | | 18 | | |
| | | | 30 | | 31 | | 24 | | | |
| 177 | | | 15 | | 23 | | 26 | | | |
| | | 22 | | 48 | | 2 | | 13 | | |
| | | 5 | | 35 | | 37 | | 49 | | |
| | 20 | | 52 | | 32 | | 50 | | 21 | |
| | 53 | | 9 | | 7 | | 10 | | 4 | |
| 31 | | 38 | | 42 | | 41 | | 43 | | 18 |
| 16 | | 36 | | 40 | | 39 | | 33 | | 28 |
| | 3 | | 12 | | 8 | | 11 | | 51 | |
| | 25 | | 54 | | 29 | | 46 | | 30 | |
| | | 47 | | 34 | | 44 | | 6 | | |
| | | 14 | | 1 | | 45 | | 19 | | |
| | | | 27 | | 24 | | 17 | | | |
| 176 | | | 17 | | 23 | | 25 | | | |
| | | 21 | | 50 | | 9 | | 14 | | |
| | | 11 | | 40 | | 42 | | 47 | | |
| | 13 | | 37 | | 12 | | 39 | | 15 | |
| | 48 | | 32 | | 2 | | 27 | | 10 | |
| 24 | | 35 | | 53 | | 54 | | 38 | | 16 |
| 22 | | 44 | | 7 | | 8 | | 43 | | 18 |
| | 3 | | 5 | | 52 | | 6 | | 51 | |
| | 30 | | 45 | | 33 | | 41 | | 31 | |
| | | 49 | | 34 | | 36 | | 4 | | |
| | | 28 | | 1 | | 46 | | 20 | | |
| | | | 19 | | 26 | | 29 | | | |
| 175 | | | 15 | | 35 | | 17 | | | |
| | | 21 | | 48 | | 1 | | 18 | | |
| | | 2 | | 40 | | 38 | | 50 | | |
| | 31 | | 49 | | 13 | | 51 | | 14 | |
| | 39 | | 12 | | 9 | | 11 | | 6 | |
| 30 | | 42 | | 52 | | 53 | | 43 | | 20 |
| 25 | | 32 | | 3 | | 4 | | 37 | | 22 |
| | 7 | | 34 | | 54 | | 27 | | 47 | |
| | 24 | | 33 | | 10 | | 36 | | 23 | |
| | | 45 | | 41 | | 44 | | 5 | | |
| | | 29 | | 8 | | 46 | | 28 | | |
| | | | 19 | | 26 | | 16 | | | |
| 174 | | | 33 | | 20 | | 16 | | | |
| | | 14 | | 3 | | 51 | | 18 | | |
| | | 30 | | 46 | | 42 | | 8 | | |
| | 26 | | 48 | | 12 | | 39 | | 15 | |
| | 7 | | 13 | | 1 | | 27 | | 47 | |
| 17 | | 50 | | 54 | 174 | 53 | | 38 | | 23 |
| 173 | | | 20 | | 12 | | 25 | | | |
| | | 18 | | 38 | | 45 | | 26 | | |
| | | 47 | | 40 | | 8 | | 27 | | |
| | 17 | | 10 | | 30 | | 42 | | 16 | |
| | 22 | | 36 | | 54 | | 34 | | 50 | |
| 35 | | 41 | | 3 | | 5 | | 4 | | 14 |
| 19 | | 41 | | 5 | | 9 | | 37 | | 25 |
| | 40 | | 11 | | 52 | | 10 | | 4 | |
| | 31 | | 45 | | 32 | | 43 | | 36 | |
| | | 6 | | 29 | | 28 | | 44 | | |
| | | 21 | | 49 | | 2 | | 35 | | |
| | | | 24 | | 34 | | 22 | | | |
| 19 | | 7 | | 53 | | 52 | | 46 | | 15 |
| | 49 | | 33 | | 6 | | 32 | | 44 | |
| | 13 | | 48 | | 31 | | 1 | | 11 | |
| | | 23 | | 2 | | 51 | | 39 | | |
| | | 28 | | 43 | | 37 | | 21 | | |
| | | | 29 | | 9 | | 24 | | | |
| 172 | | | 29 | | 15 | | 18 | | | |
| | | 8 | | 1 | | 53 | | 3 | | |
| | | 45 | | 43 | | 49 | | 7 | | |
| | 12 | | 46 | | 11 | | 42 | | 14 | |
| | 4 | | 24 | | 9 | | 25 | | 40 | |
| 48 | | 41 | | 39 | | 36 | | 44 | | 47 |
| 23 | | 5 | | 37 | | 38 | | 2 | | 17 |
| | 51 | | 26 | | 13 | | 27 | | 22 | |
| | 30 | | 10 | | 34 | | 54 | | 35 | |
| | | 50 | | 52 | | 6 | | 32 | | |
| | | 19 | | 21 | | 31 | | 33 | | |
| | | | 20 | | 28 | | 16 | | | |
| 171 | | | 41 | | 36 | | 33 | | | |
| | | 1 | | 12 | | 8 | | 15 | | |
| | | 25 | | 43 | | 46 | | 27 | | |
| | 18 | | 49 | | 26 | | 42 | | 22 | |
| | 21 | | 6 | | 10 | | 9 | | 24 | |
| 30 | | 52 | | 37 | | 38 | | 47 | | 20 |
| 16 | | 2 | | 39 | | 40 | | 3 | | 32 |
| | 50 | | 35 | | 7 | | 34 | | 45 | |
| | 19 | | 54 | | 23 | | 14 | | 31 | |
| | | 11 | | 13 | | 53 | | 44 | | |
| | | 17 | | 48 | | 5 | | 51 | | |
| | | | 28 | | 29 | | 4 | | | |
| 170 | | | 31 | | 15 | | 10 | | | |
| | | 20 | | 50 | | 51 | | 22 | | |
| | | 19 | | 9 | | 7 | | 36 | | |
| | 16 | | 41 | | 38 | | 44 | | 14 | |
| | 24 | | 27 | | 1 | | 28 | | 42 | |
| 23 | | 43 | | 54 | | 52 | | 6 | | 3 |
| 21 | | 11 | | 2 | | 8 | | 46 | | 34 |
| | 48 | | 33 | | 53 | | 30 | | 39 | |
| | 17 | | 49 | | 29 | | 5 | | 13 | |
| | | 12 | | 4 | | 45 | | 37 | | |
| | | 47 | | 32 | | 25 | | 40 | | |
| | | | 26 | | 35 | | 18 | | | |
| 169 | | | 44 | | 18 | | 19 | | | |
| | | 8 | | 2 | | 40 | | 17 | | |
| | | 42 | | 39 | | 41 | | 20 | | |
| | 24 | | 34 | | 29 | | 32 | | 23 | |
| | 4 | | 14 | | 1 | | 13 | | 31 | |
| 35 | | 51 | | 52 | | 53 | | 50 | | 33 |
| 22 | | 12 | | 3 | | 6 | | 9 | | 25 |
| | 45 | | 37 | | 54 | | 38 | | 21 | |
| | 16 | | 11 | | 15 | | 46 | | 28 | |
| | | 48 | | 49 | | 10 | | 27 | | |
| | | 26 | | 5 | | 43 | | 36 | | |
| | | | 30 | | 47 | | 7 | | | |
| 168 | | | 17 | | 22 | | 33 | | | |
| | | 25 | | 36 | | 11 | | 29 | | |
| | | 34 | | 6 | | 49 | | 42 | | |
| | 24 | | 50 | | 44 | | 4 | | 20 | |
| | 39 | | 13 | | 54 | | 14 | | 37 | |
| 18 | | 8 | | 1 | | 3 | | 51 | | 21 |
| 32 | | 48 | | 53 | | 52 | | 2 | | 30 |
| | 23 | | 45 | | 5 | | 46 | | 27 | |
| | 10 | | 7 | | 15 | | 41 | | 40 | |
| | | 35 | | 43 | | 9 | | 12 | | |
| | | 19 | | 38 | | 47 | | 28 | | |
| | | | 26 | | 16 | | 31 | | | |
| 167 | | | 21 | | 18 | | 24 | | | |
| | | 22 | | 34 | | 8 | | 31 | | |
| | | 36 | | 11 | | 51 | | 50 | | |
| | 20 | | 43 | | 45 | | 3 | | 25 | |
| | 7 | | 13 | | 1 | | 14 | | 33 | |
| 32 | | 48 | | 54 | | 53 | | 42 | | 16 |
| 27 | | 4 | | 2 | | 5 | | 6 | | 29 |
| | 49 | | 46 | | 52 | | 47 | | 41 | |
| | 19 | | 40 | | 15 | | 10 | | 28 | |
| | | 9 | | 12 | | 38 | | 35 | | |
| | | 44 | | 39 | | 37 | | 30 | | |
| | | | 23 | | 26 | | 17 | | | |
| 166 | | | 30 | | 18 | | 42 | | | |
| | | 24 | | 40 | | 3 | | 25 | | |
| | | 20 | | 45 | | 41 | | 44 | | |
| | 21 | | 7 | | 19 | | 11 | | 31 | |
| | 28 | | 39 | | 4 | | 38 | | 8 | |
| 35 | | 51 | | 52 | | 53 | | 34 | | 6 |
| 26 | | 14 | | 1 | | 2 | | 29 | | 46 |
| | 12 | | 9 | | 54 | | 10 | | 43 | |
| | 33 | | 50 | | 37 | | 47 | | 32 | |
| | | 48 | | 15 | | 16 | | 5 | | |
| | | 17 | | 13 | | 49 | | 22 | | |
| | | | 23 | | 36 | | 27 | | | |
| 165 | | | 10 | | 26 | | 38 | | | |
| | | 28 | | 39 | | 15 | | 27 | | |
| | | 23 | | 21 | | 22 | | 20 | | |
| | 36 | | 44 | | 42 | | 43 | | 41 | |
| | 9 | | 5 | | 1 | | 4 | | 8 | |
| 25 | | 48 | | 52 | 165 | 53 | | 49 | | 18 |
| 30 | | 7 | | 3 | | 2 | | 6 | | 37 |
| | 46 | | 50 | | 54 | | 51 | | 47 | |
| | 19 | | 11 | | 13 | | 12 | | 14 | |
| | | 32 | | 34 | | 33 | | 35 | | |
| | | 31 | | 40 | | 16 | | 24 | | |
| | | | 17 | | 29 | | 45 | | | |

### 4.1.1 2 by 2 triangular-shaped HTP

| | | | | |
|---|---|---|---|---|
| 50 | | 4 | | |
| | 2 | | 8 | |
| | 12 | | 11 | |
| 1 | | 13 | 50 | 3 |
| 5 | | 10 | | 7 |
| | 9 | | 6 | |
| 49 | | 4 | | |
| | 2 | | 8 | |
| | 13 | | 12 | |
| 1 | | 10 | | 3 |
| 5 | | 11 | | 7 |
| | 9 | | 6 | |
| 48 | | 4 | | |
| | 3 | | 6 | |
| | 12 | | 13 | |
| 2 | | 10 | | 1 |
| 5 | | 8 | | 9 |
| | 11 | | 7 | |
| 47 | | 4 | | |
| | 2 | | 9 | |
| | 13 | | 12 | |
| 3 | | 7 | | 1 |
| 5 | | 11 | | 10 |
| | 8 | | 6 | |
| 46 | | 4 | | |
| | 1 | | 9 | |
| | 12 | | 13 | |
| 3 | | 7 | | 2 |
| 6 | | 8 | | 5 |
| | 10 | | 11 | |
| 45 | | 7 | | |
| | 8 | | 1 | |
| | 13 | | 12 | |
| 3 | | 4 | | 10 |
| 5 | | 11 | | 2 |
| | 9 | | 6 | |
| 44 | | 5 | | |
| | 1 | | 9 | |
| | 12 | | 13 | |
| 3 | | 4 | | 2 |
| 7 | | 8 | | 6 |
| | 10 | | 11 | |
| 43 | | 5 | | |
| | 8 | | 4 | |
| | 13 | | 12 | |
| 2 | | 1 | | 10 |
| 7 | | 11 | | 3 |
| | 9 | | 6 | |
| 42 | | 3 | | |
| | 4 | | 9 | |
| | 12 | | 13 | |
| 5 | | 1 | | 2 |
| 6 | | 8 | | 7 |
| | 10 | | 11 | |
| 41 | | 9 | | |
| | 6 | | 10 | |
| | 1 | | 2 | |
| 12 | | 13 | | 4 |
| 7 | | 3 | | 11 |
| | 5 | | 8 | |
| 40 | | 9 | | |
| | 13 | | 5 | |
| | 2 | | 1 | |
| 11 | | 10 | | 12 |
| 7 | | 6 | | 8 |
| | 4 | | 3 | |
| 39 | | 7 | | |
| | 6 | | 13 | |
| | 1 | | 2 | |
| 11 | | 10 | | 4 |
| 9 | | 3 | | 12 |
| | 5 | | 8 | |
| 38 | | 10 | | |
| | 13 | | 5 | |
| | 2 | | 1 | |
| 11 | | 7 | | 12 |
| 8 | | 6 | | 9 |
| | 4 | | 3 | |
| 37 | | 10 | | |
| | 12 | | 5 | |
| | 1 | | 2 | |
| 11 | | 7 | | 13 |
| 9 | | 3 | | 4 |
| | 6 | | 8 | |
| 36 | | 10 | | |
| | 11 | | 8 | |
| | 2 | | 1 | |
| 12 | | 4 | | 13 |
| 9 | | 6 | | 5 |
| | 3 | | 7 | |
| 35 | | 10 | | |
| | 12 | | 6 | |
| | 1 | | 2 | |
| 13 | | 4 | | 11 |
| 9 | | 3 | | 7 |
| | 5 | | 8 | |
| 34 | | 10 | | |
| | 12 | | 6 | |
| | 2 | | 3 | |
| 13 | | 1 | | 11 |
| 9 | | 4 | | 7 |
| | 5 | | 8 | |

### 4.1.2 3 by 3 triangular-shaped HTP

| | | | | | | | | |
|---|---|---|---|---|---|---|---|---|
| | 121 | | | 12 | | | | |
| | | | 17 | | 10 | | | |
| | | | 24 | | 26 | | | |
| | | 4 | | 32 | | 2 | | |
| | | 23 | | 20 | | 19 | | |
| | 6 | | 18 | | 22 | | 3 | |
| | 27 | | 16 | | 15 | | 29 | |
| 21 | | 31 | | 30 | | 33 | | 7 |
| 8 | | 25 | | 14 | | 28 | | 11 |
| | 9 | | 5 | | 1 | | 13 | |
| | 120 | | | 7 | | | | |
| | | | 14 | | 17 | | | |
| | | | 27 | | 22 | | | |
| | | 1 | | 33 | | 6 | | |
| | | 19 | | 28 | | 20 | | |
| | 5 | | 12 | | 11 | | 3 | |
| | 23 | | 29 | | 30 | | 25 | |
| 8 | | 32 | | 10 | | 31 | | 9 |
| 15 | | 24 | | 21 | | 26 | | 13 |
| | 18 | | 4 | | 2 | | 16 | |
| | 119 | | | 15 | | | | |
| | | | 3 | | 13 | | | |
| | | | 26 | | 30 | | | |
| | | 4 | | 32 | | 8 | | |
| | | 19 | | 11 | | 10 | | |
| | 6 | | 27 | | 28 | | 9 | |
| | 24 | | 12 | | 16 | | 23 | |
| 20 | | 31 | | 25 | | 33 | | 22 |
| 14 | | 29 | | 17 | | 21 | | 18 |
| | 1 | | 5 | | 7 | | 2 | |
| | 118 | | | 3 | | | | |
| | | | 16 | | 20 | | | |
| | | | 23 | | 27 | | | |
| | | 8 | | 29 | | 5 | | |
| | | 10 | | 18 | | 28 | | |
| | 7 | | 30 | | 11 | | 9 | |
| | 26 | | 14 | | 13 | | 24 | |
| 2 | | 31 | | 32 | | 33 | | 1 |
| 15 | | 25 | | 12 | | 22 | | 17 |
| | 19 | | 4 | | 6 | | 21 | |
| | 117 | | | 30 | | | | |
| | | | 12 | | 8 | | | |
| | | | 16 | | 20 | | | |
| | | 5 | | 31 | | 1 | | |
| | | 25 | | 27 | | 23 | | |
| | 6 | | 13 | | 15 | | 2 | |
| | 18 | | 22 | | 26 | | 19 | |
| 7 | | 33 | | 14 | | 32 | | 9 |
| 28 | | 21 | | 24 | | 17 | | 29 |
| | 10 | | 3 | | 4 | | 11 | |
| | 116 | | | 31 | | | | |
| | | | 8 | | 9 | | | |
| | | | 20 | | 21 | | | |
| | | 5 | | 27 | | 6 | | |
| | | 29 | | 12 | | 22 | | |
| | 3 | | 23 | | 28 | | 2 | |
| | 17 | | 14 | | 15 | | 16 | |
| 26 | | 30 | | 24 | | 33 | | 10 |
| 13 | | 19 | | 25 | | 18 | | 32 |
| | 11 | | 4 | | 1 | | 7 | |
| | 115 | | | 24 | | | | |
| | | | 14 | | 25 | | | |
| | | | 9 | | 10 | | | |
| | | 3 | | 33 | | 2 | | |
| | | 29 | | 30 | | 28 | | |
| | 4 | | 11 | | 12 | | 6 | |
| | 13 | | 27 | | 18 | | 19 | |
| 8 | | 31 | | 17 | | 32 | | 7 |
| 21 | | 16 | | 23 | | 20 | | 15 |
| | 26 | | 1 | | 5 | | 22 | |
| | 114 | | | 22 | | | | |
| | | | 7 | | 26 | | | |
| | | | 14 | | 12 | | | |
| | | 5 | | 33 | | 1 | | |
| | | 25 | | 18 | | 30 | | |
| | 6 | | 19 | | 20 | | 4 | |
| | 15 | | 17 | | 16 | | 13 | |
| 9 | | 32 | | 24 | | 31 | | 8 |
| 27 | | 10 | | 29 | | 11 | | 23 |
| | 21 | | 2 | | 3 | | 28 | |
| | 113 | | | 31 | | | | |
| | | | 3 | | 22 | | | |
| | | | 15 | | 14 | | | |
| | | 8 | | 28 | | 6 | | |
| | | 19 | | 18 | | 24 | | |
| | 9 | | 25 | | 23 | | 7 | |
| | 16 | | 11 | | 10 | | 17 | |
| 1 | | 33 | | 26 | | 32 | | 2 |
| 30 | | 12 | | 27 | | 13 | | 29 |
| | 21 | | 4 | | 5 | | 20 | |
| | 112 | | | 30 | | | | |
| | | | 1 | | 24 | | | |
| | | | 14 | | 15 | | | |
| | | 4 | | 28 | | 5 | | |
| | | 27 | | 22 | | 26 | | |
| | 6 | | 17 | | 16 | | 7 | |
| | 12 | | 18 | | 20 | | 10 | |
| 3 | | 32 | | 19 | | 33 | | 2 |
| 29 | | 13 | | 21 | | 11 | | 31 |
| | 23 | | 9 | | 8 | | 25 | |
| | 111 | | | 22 | | | | |
| | | | 2 | | 31 | | | |
| | | | 11 | | 15 | | | |
| | | 8 | | 30 | | 5 | | |
| | | 25 | | 16 | | 26 | | |
| | 6 | | 21 | | 19 | | 7 | |
| | 14 | | 17 | | 18 | | 12 | |
| 3 | | 28 | | 20 | | 29 | | 1 |
| 24 | | 10 | | 27 | | 13 | | 23 |
| | 32 | | 9 | | 4 | | 33 | |
| | 110 | | | 33 | | | | |
| | | | 8 | | 1 | | | |
| | | | 23 | | 19 | | | |
| | | 17 | | 26 | | 16 | | |
| | | 5 | | 10 | | 11 | | |
| | 18 | | 29 | | 28 | | 13 | |
| | 24 | | 9 | | 7 | | 21 | |
| 2 | | 25 | | 27 | | 30 | | 3 |
| 6 | | 22 | | 12 | | 20 | | 4 |
| | 31 | | 15 | | 14 | | 32 | |
| | 109 | | | 26 | | | | |
| | | | 1 | | 11 | | | |
| | | | 24 | | 22 | | | |
| | | 19 | | 25 | | 18 | | |
| | | 8 | | 31 | | 9 | | |
| | 15 | | 2 | | 4 | | 13 | |
| | 23 | | 32 | | 33 | | 20 | |
| 3 | | 29 | | 7 | | 30 | | 5 |
| 12 | | 14 | | 6 | | 17 | | 10 |
| | 28 | | 21 | | 16 | | 27 | |
| | 108 | | | 4 | | | | |
| | | | 11 | | 27 | | | |
| | | | 19 | | 22 | | | |
| | | 24 | | 25 | | 15 | | |
| | | 1 | | 7 | | 8 | | |
| | 23 | | 32 | | 31 | | 20 | |
| | 17 | | 6 | | 2 | | 14 | |
| 10 | | 29 | | 30 | | 33 | | 12 |
| 5 | | 21 | | 9 | | 18 | | 28 |
| | 26 | | 13 | | 16 | | 3 | |

| | | | | | | | | | | | | | | | | | |
|---|---|---|---|---|---|---|---|---|---|---|---|---|---|---|---|---|---|
| | 107 | | | 10 | | | | | | 106 | | | 3 | | | | |
| | | | 1 | | 31 | | | | | | | 10 | | 30 | | | |
| | | | 23 | | 16 | | | | | | | 17 | | 21 | | | |
| | | 14 | | 26 | | 24 | | | | | 18 | | 25 | | 14 | | |
| | | 11 | | 28 | | 7 | | | | | 9 | | 33 | | 8 | | |
| | 18 | | 5 | | 6 | | 20 | | | 20 | | 4 | | 5 | | 24 | |
| | 21 | | 27 | | 29 | | 15 | | | 19 | | 28 | | 29 | | 13 | |
| 32 | | 25 | | 12 | | 30 | | 33 | 32 | | 26 | | 7 | | 27 | | 31 |
| 9 | | 17 | | 4 | | 19 | | 8 | 12 | | 16 | | 6 | | 22 | | 11 |
| | 3 | | 22 | | 13 | | 2 | | | 1 | | 23 | | 15 | | 2 | |
| | 105 | | | 33 | | | | | | 104 | | | 11 | | | | |
| | | | 1 | | 11 | | | | | | | 3 | | 33 | | | |
| | | | 16 | | 17 | | | | | | | 16 | | 14 | | | |
| | | 24 | | 27 | | 18 | | | | | 23 | | 27 | | 22 | | |
| | | 4 | | 6 | | 8 | | | | | 8 | | 29 | | 5 | | |
| | 23 | | 28 | | 29 | | 22 | | | 18 | | 1 | | 7 | | 24 | |
| | 19 | | 5 | | 7 | | 14 | | | 21 | | 30 | | 28 | | 15 | |
| 3 | | 26 | | 30 | | 25 | | 2 | 2 | | 26 | | 9 | | 25 | | 4 |
| 10 | | 15 | | 9 | | 21 | | 12 | 10 | | 13 | | 6 | | 17 | | 12 |
| | 32 | | 20 | | 13 | | 31 | | | 32 | | 20 | | 19 | | 31 | |
| | 103 | | | 18 | | | | | | 102 | | | 15 | | | | |
| | | | 29 | | 3 | | | | | | | 20 | | 16 | | | |
| | | | 24 | | 23 | | | | | | | 23 | | 25 | | | |
| | | 12 | | 6 | | 7 | | | | | 11 | | 3 | | 9 | | |
| | | 28 | | 13 | | 33 | | | | | 31 | | 4 | | 32 | | |
| | 9 | | 20 | | 21 | | 8 | | | 12 | | 30 | | 29 | | 8 | |
| | 25 | | 16 | | 14 | | 26 | | | 22 | | 5 | | 6 | | 26 | |
| 32 | | 5 | | 19 | | 1 | | 4 | 13 | | 2 | | 28 | | 1 | | 14 |
| 17 | | 22 | | 31 | | 27 | | 15 | 17 | | 27 | | 33 | | 24 | | 18 |
| | 2 | | 10 | | 11 | | 30 | | | 21 | | 7 | | 10 | | 19 | |

6.0.0 2 by 2 hexagram-shaped HTP

| | | | | | | | | |
|---|---|---|---|---|---|---|---|---|
| 129 | | | | 26 | | | | |
| | | | 17 | | 2 | | | |
| | 27 | | 42 | | 41 | | 40 | |
| 16 | | 11 | | 1 | | 3 | | 18 |
| 29 | | 32 | | 38 | | 39 | | 25 |
| | 14 | | 5 | | 7 | | 4 | |
| | 28 | | 37 | | 36 | | 31 | |
| 15 | | 13 | | 6 | | 12 | | 19 |
| 22 | | 30 | | 33 | | 34 | | 24 |
| | 21 | | 10 | | 8 | | 9 | |
| | | | 23 | | 35 | | | |
| | | | | 20 | | | | |